\documentclass[aos,MSNbibl,nameyear,dvips]{arximspdf}
\usepackage{graphicx}
\doi{10.1214/13-AOS1175A} %kopijuoti is PTS
\referstodoi{10.1214/13-AOS1175}
\volume{42}
\issue{2}
\pubyear{2014}
\firstpage{469}
\lastpage{477}

\makeatletter
\newcommand{\AR}{\mathit{AR}}
\makeatother

\begin{document}
\begin{frontmatter}

\title{Discussion: ``A significance test for the lasso''}
\runtitle{Discussion}

\begin{aug}
\author{\fnms{Peter} \snm{B\"uhlmann}\corref{}\ead[label=e1]{buhlmann@stat.math.ethz.ch}},
\author{\fnms{Lukas} \snm{Meier}\ead[label=e2]{meier@stat.math.ethz.ch}}
\and
\author{\fnms{Sara} \snm{van de Geer}\ead[label=e3]{geer@stat.math.ethz.ch}}
\runauthor{P. B\"uhlmann, L. Meier and S. van de Geer}
\affiliation{ETH Z\"urich}
\address{Seminar for Statistics\\
ETH Z\"urich\\
CH-8092 Z\"urich\\
Switzerland\\
\printead{e1}\\
\phantom{E-mail:\ }\printead*{e2}\\
\phantom{E-mail:\ }\printead*{e3}} %adresu isvedimo komanda gale!
\pdftitle{Discussion of ``A significance test for the lasso''}
\end{aug}

% HISTORY:
\received{\smonth{12} \syear{2013}}

% ABSTRACT

% KEYWORDS
% Pirmas kwd is didziosios raides
%

\end{frontmatter}

We congratulate Richard Lockhart, Jonathan Taylor, Ryan Tibshirani and
Robert Tibshirani for a thought
provoking and interesting paper on the important topic of hypothesis
testing in potentially high-dimensional settings.

%s1 #&#
\section{A short description of the test procedure}
We start by presenting the proposed test procedure in a slightly different
form than in the paper. Let
\[
\hat\beta(\lambda):= \arg\min \tfrac{1}{2} \| y - X\beta
\|_2^2 + \lambda\| \beta\|_1
\]
be the Lasso estimator with tuning parameter equal to $\lambda$.
The paper uses the Lasso path $\{ \hat\beta(\lambda)\dvtx  \lambda> 0
\} $ to construct a test statistic for the significance of certain
predictor variables.

For a subset $S \subseteq\{ 1, \ldots, p \}$, let $\hat\beta_S
(\lambda)$ be the Lasso
solution using only the variables in $S$:
\[
\hat\beta_S (\lambda):= \arg\min_{\beta_S \in\mathbb{R}^{|S|}
}
{1 \over2} \| y - X_S\beta_S
\|_2^2 + \lambda\| \beta_S \|_1.
\]
The covariance test is based on the difference
\begin{eqnarray*}
T ( S, \lambda) &:=& \bigl[ \bigl\| y - X_S \hat\beta_S (
\lambda) \bigr\|_2^2 + \lambda\bigl\| \hat\beta_S (
\lambda) \bigr\|_1 \bigr] / \sigma^2
\\
&&{}- \bigl[ \bigl\| y - X \hat
\beta(\lambda) \bigr\|_2^2 + \lambda\bigl\| \hat\beta(\lambda)
\bigr\|_1 \bigr]/\sigma^2.
\end{eqnarray*}
If $T(S, \lambda)$ is large, then the solution using only the values
in $S$ does not have a very good
fit, and this may support evidence against the hypothesis $H_S\dvtx  A^*
\subseteq S$, where
$A^* = \operatorname{support} ( \beta^*)$ is the true active set.\vspace*{1pt}

Let $\infty=:\hat\lambda_0 > \hat\lambda_1 \ge\hat\lambda_2 \ge \cdots$ be the knots
of $\hat\beta(\lambda)$. For $k \ge1$, let $\widehat A_k:= \operatorname{support}
( \hat\beta(\hat\lambda_{k} ))$.
We put ``hats'' on these quantities to stress that they are random
variables depending (only) on the data.\vspace*{1pt}
Thus,
$T(\widehat A_k, \hat\lambda_k ) =0 $ and by continuity arguments also
$T(\widehat A_{k-1}, \hat\lambda_k ) = 0 $. The authors suggest to use
the test statistic
\[
T_k:= T( \widehat A_{k-1}, \hat\lambda_{k+1} )
\]
for\vspace*{1pt} the hypotheses $H_{\widehat A_{k-1} }$. They derive the interesting result
that under certain conditions, the test statistic has an asymptotic
exponential distribution.

%
%f1 #&#
\begin{figure}[b]

\includegraphics{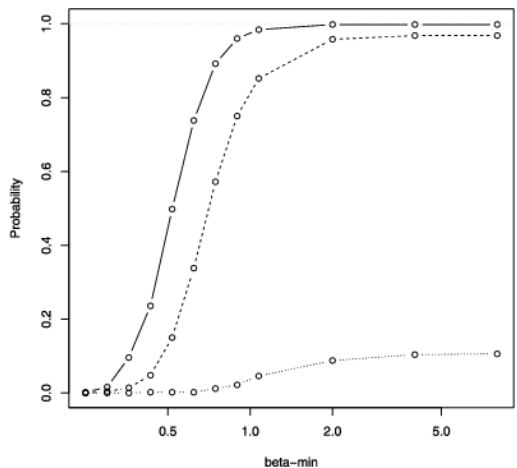}

\caption{Empirical probabilities (500 simulation runs) for the event that
all truly active coefficients are identified in the first $k_0$ steps of
the Lasso path. A Gaussian $\AR(1)$ design matrix with $n=100$, $p = 1000$ and
$\rho= 0.5$ is used (i.e., the population covariance matrix $\Sigma$ is
Toeplitz with $\Sigma_{ij} = \rho^{|i-j|})$. The active set has size
$k_0=3$ (solid), $k_0=5$
(dashed) and $k_0=10$ (dotted). The active coefficients are placed at
random positions and all coefficients are of the same size
(beta-min). The error variance $\sigma^2$ is set to 1.}\label{figprobB}
\end{figure}
%

%s2 #&#
\section{A ``conditional'' test}

Fixing the value of $k$, the test is a conditional test for $H_S$ given
that $\widehat{A}_{k-1} = S$ (the event one conditions on is denoted below and
in the paper by $B$; the paper presents two versions, in Sections~3.2
and 4.2, resp.). Such kind of a test is uncommon: the usual
form of a conditional test is to condition on an observable event, for
example, when conditioning on an ancillary statistics
[cf. \citet{ghoshetal10}]. Here, however, the
conditioning event $B$, that
all active variables enter the Lasso path first, is \emph{unobserved}.

The difficulty with such an unobserved event is treated in the paper by
imposing sufficient conditions such that $\mathbb{P}[B] \to1$
asymptotically and,
therefore, one can simply ignore the effect of conditioning. The imposed
conditions are rather restrictive: in particular, they include a
``beta-min'' assumption requiring that the nonzero regression coefficients
are sufficiently large in absolute value. We illustrate in Figure~\ref{figprobB} that the lower bound for the nonzero coefficients
(beta-min) has to
be large or very large in order that the active set is correctly identified
right at the first steps of the Lasso path [the latter is the conditioning
event $B$ as in Section~3.2 of the paper while in Section~4.2 of the
paper, a slightly different version of $B$ is presented; we believe that
the quantitative differences in terms of $P(B)$ are small].
Based on this observation, we imagine that the obtained
limiting distribution in Theorems~1 and~3 often does not approximately
capture the conditional distribution of the test statistics (when
conditioning on the
event $B$), and there is no strong guarantee that the obtained $p$-values
would be approximately correct in practical settings. It would be
interesting to work out a correction factor which would take into
account that
$\mathbb{P}[B]$ is not close to one: we do not know how this could be achieved.

%s2.1 #&#
\subsection{Interpretation of the $p$-values}\label{subsecinterp}

A correct interpretation of the proposed $p$-values
from the covariance test seems not straightforward. First, these $p$-values
are \emph{not} justified to lead to significance statements for fixed
variables (or hypotheses) since the test is a conditional test. For
example, for the wine data in the right
column of Table~5 in the paper, the $p$-value for the variable
``\texttt{pH}'' should not be interpreted in the classical sense based
on a
fixed null-hypothesis\vspace*{1pt} $\beta^*_{\mathrm{pH}} = 0$. In
many scientific applications and contexts, such classical $p$-values are
desired (maybe after adjustment for multiple testing), and we think that
the covariance test does not really provide such $p$-values; in fact,
the authors never make such a claim. Reasons for the statement above
include: (i) the covariance test only assigns significance of ``the $k$th
variable entering the Lasso path,'' but since the $k$th variable is
random (possibly even when $\mathbb{P}[B] \to1$), there seems to be
an issue to map
the $k$th variable to a fixed variable, such as ``\texttt{pH}'' or
``\texttt{alcohol}''; (ii) in view that $\mathbb{P}[B]$ might be far
away from
one as illustrated in Figure~\ref{figprobB}, the interpretation
should be
conditional, that is, ``given that
all active variables enter the Lasso path first''; and such a conditional
interpretation of a $p$-value seems somewhat awkward. We briefly outline
in Section~\ref{secaltern} some alternative methods which are
mathematically justified for classical (fixed hypotheses) $p$-values in a
high-dimensional context.

Our question to the authors is how to interpret the $p$-values in
practice. In view of available software, there is a substantial risk that
practitioners blindly use and interpret the obtained $p$-values as usual (for
fixed hypotheses), and
hence, some guidance for proper use and interpretation would be very useful.

%s3 #&#
\section{The assumptions}

The authors require a condition on the design matrix and a beta-min
assumption. These assumptions are used to guarantee that the conditioning
event $B$, namely that the first $k_0$ variables entering the Lasso path
contain the active set, has large probability.

In Theorem~3 of the paper, an irrepresentable condition [\citet{zhaoyu06}] is
assumed. Let
$A_0 \supseteq A^*$ and let
\[
\eta> \max_{j \notin A_0} \sup_{\| \tau_{A_0} \|_{\infty} \le1 } \bigl|
X_j^T X_{A_0} \bigl( X_{A_0}^T
X_{A_0} \bigr)^{-1} \tau_{A_0} \bigr|.
\]
We assume the irrepresentable condition $\eta\le1$.
From Exercise 7.5 in \citet{pbvdg11} we know that for
$ \lambda_{\eta}:= \lambda_{\epsilon} (1+ \eta) / (1 - \eta) $
we have \mbox{$\widehat A ( \lambda_{\eta}) \subseteq A_0$.} Here,
\[
\lambda_{\epsilon} = \max_{1 \le j \le p} \bigl| \langle\epsilon,
X_j \rangle\bigr|.
\]
Define now
\[
\hat k_{\eta}:= \max\{ k\dvtx  \lambda_k \ge
\lambda_{\eta} \}.
\]
Thus, with large probability,
\[
A_* \subseteq\widehat A_{\hat k_{\eta} } \subseteq A_0.
\]
We imagine moreover that in practice one would follow the
Lasso path and steps as soon as the test accepts $\widehat A_{k-1}$.
Define therefore $\hat k$ as being the first
$k$ for which the hypothesis $H_{\widehat A_{k-1}}$ is accepted.
In line with the paper, one then assumes $\widehat A_{\hat k-1}\supseteq A_0$,
and then with probability approximately $1- \alpha$, $\widehat A_{\hat
k-1} = A_0$.
Alternatively, applying this argument to $\widehat A_{\hat k_{\eta} }$
(which is allowed since $A_* \subseteq\widehat A_{\hat k_{\eta} }$) we get
$\hat k_{\eta} \ge\hat k-1$ with probability approximately $1-\alpha
$ and then
we end up with $A_* \subseteq\widehat A_{\hat k-1} = \widehat A_{\hat k_{\eta
}} = A_0 $.
%Then
%$$ \widehat A_{\hat k_{\eta}} \subseteq A_0 \subseteq\widehat A_{\hat k -1}
%.$$
%If the nonzero coefficients of $\beta^*$ are sufficiently larger than
%$\lambda_{\eta}$
%one has in fact $\widehat A_{\hat k_0} \supseteq A^*$. And with probability
%approximately $1- \alpha$
%
%???SARA: WHAT IS $\alpha$; SHOULD WE WRITE And with high
%probability????
%we have $\hat k=1$.
%
%???SARA: SHOULD IT BE $\hat{k} = k_0$???
%
%So we get that with large probability
%$$A^* \subseteq\widehat A_{\hat k_{\eta}} \subseteq A_0 = \widehat A_{\hat k
%-1}.$$

A related screening property of the Lasso is known
[\citet{pbvdg11}, cf. Chapter~2.5]: for $\lambda
\asymp\sqrt{\log(p)/n}$,
%e1 #&#
\begin{equation}
\label{screening} \mathbb{P}\bigl[\widehat{A}(\lambda) \supseteq A^*\bigr]
\to1,
\end{equation}
assuming the compatibility condition on the design and a beta-min
assumption. We note that the compatibility condition is weaker than the
irrepresentable condition mentioned above [\citet{van2009conditions}].

The authors argue in their Remark 4 that the beta-min assumption can be
relaxed. Such kind of
a relaxation is given in \citet{pbmand13}, \mbox{assuming} a zonal
assumption allowing that some but not too many nonzero regression coefficients
are small. It is also shown that zonal assumptions are
necessary for validity of a sampling splitting procedure [\citet{WR08}], and
we believe that a justification of the covariance test also necessarily
needs some version of zonal assumptions. We remark that ``in practice,''
achieving a
statement as in (\ref{screening}) or saying that $\mathbb{P}[B]
\approx1$ (as in
the paper) seems often unrealistic, as illustrated in Figure~\ref{figprobB} and in \citet{pbmand13}.

%s3.1 #&#
\subsection{Hypothesis testing and assumptions on \texorpdfstring{$\beta^*$}{beta*}}

In view of the fact that assumptions about $\beta^*$ are (have to be)
made, the
covariance test is exposed to the following somewhat undesirable issue. A
significance test should
\emph{find out} whether a regression coefficient is sufficiently
large. Thus, a
zonal or beta-min assumption rules out the \emph{essence of the
question} by
assuming that
most or all nonzero coefficients are large. We note that (multi) sample
splitting techniques [\citet{WR08}, \citet{memepb09}] for hypothesis testing in
high-dimensional scenarios suffer from the same problem. The procedure
outlined in Section~\ref{secaltern} does not make such zonal or beta-min
assumptions.

%s4 #&#
\section{The power of the covariance test}

The paper does not make any claim about the power of the test nor does it
include a comparison with other methods; regarding the latter, see Section~\ref{subseccompar}.

Under the beta-min assumption, a theoretical study of the test's power is
uninteresting: asymptotically, the power of the test is approaching
one. Nontrivial power statement require that the nonzero regression
coefficients are in the $1/\sqrt{n}$ range but this is excluded by the
imposed beta-min assumption.

The following thoughts might lead to some insights for which
scenarios the covariance test is expected to perform (reasonably) well. In
an alternative and simplified setup, one could think of using a refitting
procedure to test significance.
Let
\[
\hat\beta_S:= \hat\beta_S (0) = \arg\min
_{\beta_S \in\mathbb
{R}^{|S|} } \| y - X_S \beta_S
\|_2^2
\]
and for $\widetilde S \supseteq S $
\[
T( S, \widetilde S):= \| y - X_S \hat\beta_S
\|_2^2/ \sigma^2 - \| y - X_{\widetilde S}
\hat\beta_{\widetilde S} \|_2^2/ \sigma^2
\]
\[
= \bigl( \langle y, X_{\widetilde S} \hat\beta_{\widetilde S} \rangle- \langle
y, X_{ S} \hat\beta_{S} \rangle\bigr) /
\sigma^2.
\]
An alternative test statistic would then be $T( \widehat A_{k-1}, \widehat A_{k} )
$. In the
case of orthogonal design, we get
%e2 #&#
\begin{equation}
\label{add1} T_k = \bigl( \hat{\lambda}_k^2
- \hat{\lambda}_k \hat{\lambda}_{k+1} \bigr) /
\sigma^2,\qquad T( \widehat A_{k-1}, \widehat A_{k} ) = \hat{
\lambda}_{k}^2 /\sigma^2.
\end{equation}
Obviously,\vspace*{1pt} if we fix $S$ and $j \notin S $, we get
$T( S, S \cup\{ j \} ) = (\langle y, X_j \rangle)^2 / \sigma^2 $
which has
under $H_S$ a $\chi^2 (1)$
distribution. If $\mathbb{P}( \widehat A_{k-1} \supseteq A^*) \rightarrow
1$, then for
each $j \notin\widehat A_{k-1}$, $T( \widehat A_{k-1}, \widehat A_{k-1} \cup\{ j
\} )$ is asymptotically
$\chi^2 (1)$. However, $T_k$ and $T( \widehat A_{k-1}, \widehat A_{k})$ are
tests where the decision which variable is to be tested for significance
depends on the data.
For the case of orthogonal design and $A^* = \varnothing$, we have
$\mathbb{P}_{H_{\varnothing}} ( \widehat A_{k-1} \supseteq A^*) = 1$, and
$T( \widehat A_{k-1}, \widehat A_{k})$ is approximately
distributed as the $k$th order statistic of a sample from
a $\chi^2 (1)$-distribution in decreasing order. For $k=1$ (say), the
statistic $T_1$ has a different scaling
under the hypothesis $H_{\varnothing}\dvtx  A^* = \varnothing$
because the order statistics behave like
\[
T( \varnothing, \widehat A_{1} ) ={\mathcal O}_{\mathbb{P}_{H_{\varnothing
}}} (\log n )
\]
($p =n$ in the orthonormal case) whereas $T_1$ has asymptotically an
exponential distribution, a nice fact proved in the paper, so that $T_1 =
{\mathcal O}_{\mathbb{P}_{H_{\varnothing}}} (1) $. This means
that $T_1$ has more power to detect alternatives of the form $H_{\{ j \}
}\dvtx  A^* = \{ j \} $. But it may
have less power for alternatives of the form $H_{\{ j_1, j_2 \}}\dvtx  A^* =
\{ j_1, j_2 \} $. Under this
alternative, $A^* \neq \widehat A_0 $ and if the two nonzero coefficients
are very
close together it will downscale the statistic $T_1$. This can also be
seen from the expression (\ref{add1}): if the nonzero coefficients are
similar, then
\[
\hat{\lambda}_{k-1} \approx\hat{\lambda}_k,
\]
which leads to small values for $T_k$ while this has no (substantial)
effect on $T(\widehat{A}_{k-1},\widehat{A}_k)$: thus, the covariance test
might be
subideal for detection of coefficient vectors whose individual nonzero
coefficients are similar (as in the simulated examples in the
paper and in Section~\ref{subseccompar}). It would be interesting to
better understand the regimes where the covariance test has strong and weak
power.

%s5 #&#
\section{Alternative methods}\label{secaltern}

Other methods leading to $p$-values for fixed hypotheses $H_{0,j}\dvtx  \beta^*_j
= 0$ have been proposed in earlier work
[\citet{WR08}, \citet{memepb09}, \citet{minnieretal11},
\citet{pb13}, \citet{chatter13}, \citet{zhangzhang11}].
We outline
here the method from \citet{zhangzhang11} which has been further
analyzed in
\citet{vdgetal13} and \citet{jamo13b}. The idea is to
desparsify the Lasso,
resulting in a new estimator $\hat{b}$ which is not sparse. Due to
nonsparsity, this new $\hat{b}$ will not be suitable for prediction in
high-dimensional settings, but its $j$th component $\hat{b}_j$ is
asymptotically optimal for the low-dimensional target $\beta^*_j$ of
interest:
%e3 #&#
\begin{equation}
\label{pivot} \sqrt{n} \bigl(\hat{b}_j - \beta^*_j
\bigr) \Rightarrow{\mathcal N}\bigl(0,\sigma_{\varepsilon}^2
v_j\bigr),
\end{equation}
where $\sigma_{\varepsilon}^2 v_j$ is the Cramer--Rao lower bound.
Such a result
needs some assumptions on the design and sparsity of $\beta^*$ but no
further restrictions on $\beta^*$ in terms of zonal or beta-min assumptions
[\citet{vdgetal13}, \citet{jamo13b}]. Thus,
we are in the semiparametric framework,
where we can optimally estimate a low-dimensional parameter of interest in
presence of a very high-dimensional nuisance parameter $\eta=
\{\beta^*_k; k \neq j\}$: notably, we have the $1/\sqrt{n}$ convergence
rate, even when $p \gg n$, and the best possible constant in the asymptotic
variance.

The analysis in \citet{vdgetal13} also shows that (\ref{pivot}) holds
\emph{uniformly} over all sparse parameter vectors $\beta^*$ and, therefore,
the obtained confidence intervals and tests are honest. This is not the
case when using a residual-based bootstrap in \citet{chatter13} which
exhibits the unpleasant super-efficiency phenomenon. As a consequence,
post-model selection techniques [\citet{leebpoetsch03}, \citet{berketal13}] are not
necessary to construct valid, and in fact most powerful, hypothesis testing.

%s5.1 #&#
\subsection{A small empirical comparison}\label{subseccompar}

We present here some result from a small simulation study based on a
similar model as the Gaussian $\AR(1)$ model in the paper with
$\rho=0.5$. We use an
active set $A^*$ of size 10, where the active coefficients are placed
at random
positions and all have the same size. A total of 500 simulation runs are
performed for each scenario.

We consider two-sided testing of individual hypotheses $H_{0,j}\dvtx  \beta
^*_j = 0$,
possibly with adjustment for multiple testing using the Bonferroni--Holm
procedure to control the familywise error rate.

%
%t1 #&#
\begin{table}[b]
\tabcolsep=0pt
\caption{(Empirical) familywise error rate (FWER) and average number of
true positives (TP) for desparsified Lasso (\textup{de-spars}) and both
approaches of the covariance test ($\mathrm{cov}$ and $\mathrm{cov.pval}$). The different rows correspond to coefficient
size 0.5, 1, 2 and 4 (top to bottom). Sample size $n=100$ and dimension $p=80$}\label{tabsimA}
\begin{tabular*}{\tablewidth}{@{\extracolsep{\fill}}@{}lccccc@{}}
\hline
\textbf{FWER\tsub{de-spars}} & \textbf{TP\tsub{de-spars}} & \textbf{FWER\tsub{cov}} & \textbf{TP\tsub{cov}} & \textbf{FWER\tsub{cov.pval}} & \textbf{TP\tsub{cov.pval}}
\\
\hline
0.042 & 2.626 & 0.072 & 1.304 & 0.020 & 0.736 \\
0.056 & 7.104 & 0.124 & 2.884 & 0.064 & 3.770 \\
0.064 & 9.116 & 0.284 & 5.992 & 0.210 & 7.556 \\
0.064 & 9.478 & 0.426 & 8.394 & 0.298 & 9.324\\
\hline
\end{tabular*}
\end{table}

The covariance test is used in the following two ways. A first approach
(denoted by $\mathrm{cov}$) is to follow the Lasso path until the
first time the
(unadjusted) $p$-value of the covariance test is nonsignificant and
declare all
corresponding predictor variables as significant which entered before
such a
nonsignificance flag of the covariance test.
A second
approach (denoted by $\mathrm{cov.pval}$) is to assign those predictors
that remain
in the Lasso path until the end, the $p$-value of the covariance test when
they last entered the path. The $p$-values from this second approach are then
corrected for multiple testing using the Bonferroni--Holm procedure. The
second approach might be inappropriate; see
also our discussion in Section~\ref{subsecinterp} above pointing to the fact that the covariance test
does not seem to test the hypotheses $H_{0,j}$; but for the sake of
comparison (and practical use of the covariance test), we see no immediate
other way to use the covariance test for constructing $p$-values for fixed
hypotheses.
For the $p>n$ situations, we use for all methods the variance estimator
from the scaled Lasso [\citet{sunzhang11}].

The results for $n = 100$ and $p = 80$ are reported in Table~\ref{tabsimA} and the results for $n=100$ and $p=200$ can be found in
Table~\ref{tabsimB}. In all settings, the desparsified Lasso method
reliably controls the familywise error rate. In the $n>p$ setting, the
covariance test has reasonable power at the cost of no control of the
familywise
error rate. In the $p>n$ setting, the covariance test seems to be very
conservative.
%This might be due to an unsuitable estimate of the error
%variance $\sigma$.

%
%t2 #&#
\begin{table}
\tabcolsep=0pt
\caption{(Empirical) familywise error rate (FWER) and average number of
true positives (TP) for desparsified Lasso (\textup{de-spars})
and both approaches of the covariance test ($\mathrm{cov}$ and
$\mathrm{cov.pval}$). The different rows correspond to coefficient
size 0.5, 1, 2 and 4 (top to bottom). Sample size $n=100$ and dimension $p=200$}\label{tabsimB}
\begin{tabular*}{\tablewidth}{@{\extracolsep{\fill}}@{}lccccc@{}}
\hline
\textbf{FWER\tsub{de-spars}} & \textbf{TP\tsub{de-spars}} & \textbf{FWER\tsub{cov}} & \textbf{TP\tsub{cov}} & \textbf{FWER\tsub{cov.pval}} & \textbf{TP\tsub{cov.pval}}
\\
\hline
0.030 & 1.320 & 0.012 & 0.416 & 0.002 & 0.120 \\
0.046 & 3.304 & 0.010 & 0.632 & 0.004 & 0.274 \\
0.052 & 4.934 & 0.018 & 0.956 & 0.006 & 0.860 \\
0.060 & 5.594 & 0.032 & 1.550 & 0.018 & 1.884\\
\hline
\end{tabular*}
\end{table}

%s6 #&#
\section{Conclusions}

The authors present a novel and original idea of significance testing for
``random hypotheses.'' In complex data scenarios, the strategy of
considering ``data-driven'' hypotheses is certainly interesting, and the
topic deserves further attention. The proposed solution to deal with such
``random hypotheses'' is based on strong beta-min or zonal assumptions, and
this is somewhat unsatisfactory. The idea of taking the selection effect
into account appears in other work for controlling the (Bayesian) false
discovery rate [cf. \citet{beye05}, \citet{hwangzhao13}].
We think that recent alternative approaches, as outlined in Section~\ref{secaltern}, are often more powerful and simpler to interpret
when adopting
the classical framework of (multiple) fixed hypotheses testing. It is
an open
question though whether the classical framework is the most appropriate
tool for assigning ``relevance'' of single or
groups of variables in complex or high-dimensional settings.

% zodis "Acknowledgments" paliekamas pagal autoriu

%suskaldyti doi

% imsref loaded by linak, 2014-03-18 09:37:27
% imsref loaded by linak, 2014-03-18 09:37:53
% imsref loaded by linak, 2014-03-18 10:01:02
% imsref loaded by linak, 2014-03-18 10:07:26
% imsref loaded by linak, 2014-03-18 10:10:38

\printaddresses


\begin{thebibliography}{18}

%b1 #&#
\bibitem[\protect\citeauthoryear{Benjamini and Yekutieli}{2005}]{beye05}
\begin{barticle}[mr]
\bauthor{\bsnm{Benjamini},~\bfnm{Yoav}\binits{Y.}} \AND
\bauthor{\bsnm{Yekutieli},~\bfnm{Daniel}\binits{D.}}
(\byear{2005}).
\btitle{False discovery rate-adjusted multiple confidence intervals for selected parameters}.
\bjournal{J. Amer. Statist. Assoc.}
\bvolume{100}
\bpages{71--93}.%
\bid{doi={10.1198/016214504000001907}, issn={0162-1459}, mr={2156820}}%
\bptnote{check related}%
\end{barticle}%
\bptok{imsref}%
% NOT OUTPUTED:
% issn = 0162-1459
% url = http://dx.doi.org/10.1198/016214504000001907
% number = 469
% coden = JSTNAL
% fjournal = Journal of the American Statistical Association
\endbibitem

%b2 #&#
\bibitem[\protect\citeauthoryear{Berk et~al.}{2013}]{berketal13}
\begin{barticle}[mr]
\bauthor{\bsnm{Berk},~\bfnm{Richard}\binits{R.}},
\bauthor{\bsnm{Brown},~\bfnm{Lawrence}\binits{L.}},
\bauthor{\bsnm{Buja},~\bfnm{Andreas}\binits{A.}},
\bauthor{\bsnm{Zhang},~\bfnm{Kai}\binits{K.}} \AND
\bauthor{\bsnm{Zhao},~\bfnm{Linda}\binits{L.}}
(\byear{2013}).
\btitle{Valid post-selection inference}.
\bjournal{Ann. Statist.}
\bvolume{41}
\bpages{802--837}.
\bid{doi={10.1214/12-AOS1077}, issn={0090-5364}, mr={3099122}}
\end{barticle}
\bptok{imsref}%
% NOT OUTPUTED:
% issn = 0090-5364
% url = http://dx.doi.org/10.1214/12-AOS1077
% number = 2
% coden = ASTSC7
% fjournal = The Annals of Statistics
\endbibitem

%b3 #&#
\bibitem[\protect\citeauthoryear{B{\"u}hlmann}{2013}]{pb13}
\begin{barticle}[auto:STB|2014/02/12|14:17:21]
\bauthor{\bsnm{B{\"u}hlmann},~\bfnm{P.}\binits{P.}}
(\byear{2013}).
\btitle{Statistical significance in high-dimensional linear models}.
\bjournal{Bernoulli}
\bvolume{19}
\bpages{1212--1242}.
\end{barticle}
\bptok{imsref}%
\endbibitem

%b4 #&#
\bibitem[\protect\citeauthoryear{B{\"u}hlmann and Mandozzi}{2013}]{pbmand13}
\begin{bmisc}[auto:STB|2014/02/12|14:17:21]
\bauthor{\bsnm{B{\"u}hlmann},~\bfnm{P.}\binits{P.}} \AND
\bauthor{\bsnm{Mandozzi},~\bfnm{J.}\binits{J.}}
(\byear{2013}).
\bhowpublished{High-dimensional variable screening and bias in subsequent inference,
with an empirical comparison. \textit{Comput. Statist.}
DOI:\doiurl{10.1007/s00180-013-0436-3}.}
\end{bmisc}
\bptok{imsref}%
% NOT OUTPUTED:
% sortkey = Buhlmann(2013
\endbibitem

%b5 #&#
\bibitem[\protect\citeauthoryear{B{\"u}hlmann and van~de Geer}{2011}]{pbvdg11}
\begin{bbook}[mr]
\bauthor{\bsnm{B{\"u}hlmann},~\bfnm{Peter}\binits{P.}} \AND
\bauthor{\bsnm{van~de Geer},~\bfnm{Sara}\binits{S.}}
(\byear{2011}).
\btitle{Statistics for High-Dimensional Data: Methods, Theory and Applications}.
\bpublisher{Springer},
\blocation{Heidelberg}.
\bid{doi={10.1007/978-3-642-20192-9}, mr={2807761}}
\end{bbook}
\bptok{imsref}%
% NOT OUTPUTED:
% isbn = 978-3-642-20191-2
% url = http://dx.doi.org/10.1007/978-3-642-20192-9
% fpage = xviii+556
\endbibitem

%b6 #&#
\bibitem[\protect\citeauthoryear{Chatterjee and Lahiri}{2013}]{chatter13}
\begin{barticle}[mr]
\bauthor{\bsnm{Chatterjee},~\bfnm{A.}\binits{A.}} \AND
\bauthor{\bsnm{Lahiri},~\bfnm{S.~N.}\binits{S.~N.}}
(\byear{2013}).
\btitle{Rates of convergence of the adaptive {LASSO} estimators to the oracle distribution and higher order refinements by the bootstrap}.
\bjournal{Ann. Statist.}
\bvolume{41}
\bpages{1232--1259}.
\bid{doi={10.1214/13-AOS1106}, issn={0090-5364}, mr={3113809}}
\end{barticle}
\bptok{imsref}%
% NOT OUTPUTED:
% issn = 0090-5364
% url = http://dx.doi.org/10.1214/13-AOS1106
% number = 3
% fjournal = The Annals of Statistics
\endbibitem

%b7 #&#
\bibitem[\protect\citeauthoryear{Ghosh, Reid and Fraser}{2010}]{ghoshetal10}
\begin{barticle}[mr]
\bauthor{\bsnm{Ghosh},~\bfnm{M.}\binits{M.}},
\bauthor{\bsnm{Reid},~\bfnm{N.}\binits{N.}} \AND
\bauthor{\bsnm{Fraser},~\bfnm{D.~A.~S.}\binits{D.~A.~S.}}
(\byear{2010}).
\btitle{Ancillary statistics: A review}.
\bjournal{Statist. Sinica}
\bvolume{20}
\bpages{1309--1332}.
\bid{issn={1017-0405}, mr={2777327}}
\end{barticle}
\bptok{imsref}%
% NOT OUTPUTED:
% issn = 1017-0405
% number = 4
% fjournal = Statistica Sinica
\endbibitem

%b8 #&#
\bibitem[\protect\citeauthoryear{Hwang and Zhao}{2013}]{hwangzhao13}
\begin{barticle}[mr]
\bauthor{\bsnm{Hwang},~\bfnm{J.~T.~Gene}\binits{J.~T.~G.}} \AND
\bauthor{\bsnm{Zhao},~\bfnm{Zhigen}\binits{Z.}}
(\byear{2013}).
\btitle{Empirical {B}ayes confidence intervals for selected parameters in high-dimensional data}.
\bjournal{J. Amer. Statist. Assoc.}
\bvolume{108}
\bpages{607--618}.
\bid{doi={10.1080/01621459.2013.771102}, issn={0162-1459}, mr={3174645}}
\end{barticle}
\bptok{imsref}%
% NOT OUTPUTED:
% issn = 0162-1459
% url = http://dx.doi.org/10.1080/01621459.2013.771102
% number = 502
% fjournal = Journal of the American Statistical Association
\endbibitem

%b9 #&#
\bibitem[\protect\citeauthoryear{Javanmard and Montanari}{2013}]{jamo13b}
\begin{bmisc}[auto:STB|2014/02/12|14:17:21]
\bauthor{\bsnm{Javanmard},~\bfnm{A.}\binits{A.}} \AND
\bauthor{\bsnm{Montanari},~\bfnm{A.}\binits{A.}}
(\byear{2013}).
\bhowpublished{Confidence intervals and hypothesis testing for high-dimensional regression.
Preprint. Available at \arxivurl{arXiv:1306.3171}.}
\end{bmisc}
\bptok{imsref}%
% NOT OUTPUTED:
% sortkey = Javanmard(2013
\endbibitem

%b10 #&#
\bibitem[\protect\citeauthoryear{Leeb and P{\"o}tscher}{2003}]{leebpoetsch03}
\begin{barticle}[mr]
\bauthor{\bsnm{Leeb},~\bfnm{Hannes}\binits{H.}} \AND
\bauthor{\bsnm{P{\"o}tscher},~\bfnm{Benedikt~M.}\binits{B.~M.}}
(\byear{2003}).
\btitle{The finite-sample distribution of post-model-selection estimators and uniform versus nonuniform approximations}.
\bjournal{Econometric Theory}
\bvolume{19}
\bpages{100--142}.
\bid{doi={10.1017/S0266466603191050}, issn={0266-4666}, mr={1965844}}
\end{barticle}
\bptok{imsref}%
% NOT OUTPUTED:
% issn = 0266-4666
% url = http://dx.doi.org/10.1017/S0266466603191050
% number = 1
% fjournal = Econometric Theory
\endbibitem

%b11 #&#
\bibitem[\protect\citeauthoryear{Meinshausen, Meier and B{\"u}hlmann}{2009}]{memepb09}
\begin{barticle}[mr]
\bauthor{\bsnm{Meinshausen},~\bfnm{Nicolai}\binits{N.}},
\bauthor{\bsnm{Meier},~\bfnm{Lukas}\binits{L.}} \AND
\bauthor{\bsnm{B{\"u}hlmann},~\bfnm{Peter}\binits{P.}}
(\byear{2009}).
\btitle{{$p$}-values for high-dimensional regression}.
\bjournal{J. Amer. Statist. Assoc.}
\bvolume{104}
\bpages{1671--1681}.
\bid{doi={10.1198/jasa.2009.tm08647}, issn={0162-1459}, mr={2750584}}
\end{barticle}
\bptok{imsref}%
% NOT OUTPUTED:
% issn = 0162-1459
% url = http://dx.doi.org/10.1198/jasa.2009.tm08647
% number = 488
% coden = JSTNAL
% fjournal = Journal of the American Statistical Association
\endbibitem

%b12 #&#
\bibitem[\protect\citeauthoryear{Minnier, Tian and Cai}{2011}]{minnieretal11}
\begin{barticle}[mr]
\bauthor{\bsnm{Minnier},~\bfnm{Jessica}\binits{J.}},
\bauthor{\bsnm{Tian},~\bfnm{Lu}\binits{L.}} \AND
\bauthor{\bsnm{Cai},~\bfnm{Tianxi}\binits{T.}}
(\byear{2011}).
\btitle{A perturbation method for inference on regularized regression estimates}.
\bjournal{J. Amer. Statist. Assoc.}
\bvolume{106}
\bpages{1371--1382}.
\bid{doi={10.1198/jasa.2011.tm10382}, issn={0162-1459}, mr={2896842}}
\end{barticle}
\bptok{imsref}%
% NOT OUTPUTED:
% issn = 0162-1459
% url = http://dx.doi.org/10.1198/jasa.2011.tm10382
% number = 496
% coden = JSTNAL
% fjournal = Journal of the American Statistical Association
\endbibitem

%b13 #&#
\bibitem[\protect\citeauthoryear{Sun and Zhang}{2012}]{sunzhang11}
\begin{barticle}[mr]
\bauthor{\bsnm{Sun},~\bfnm{Tingni}\binits{T.}} \AND
\bauthor{\bsnm{Zhang},~\bfnm{Cun-Hui}\binits{C.-H.}}
(\byear{2012}).
\btitle{Scaled sparse linear regression}.
\bjournal{Biometrika}
\bvolume{99}
\bpages{879--898}.
\bid{doi={10.1093/biomet/ass043}, issn={0006-3444}, mr={2999166}}
\end{barticle}
\bptok{imsref}%
% NOT OUTPUTED:
% issn = 0006-3444
% url = http://dx.doi.org/10.1093/biomet/ass043
% number = 4
% fjournal = Biometrika
\endbibitem

%b14 #&#
\bibitem[\protect\citeauthoryear{van~de Geer and B{\"u}hlmann}{2009}]{van2009conditions}
\begin{barticle}[mr]
\bauthor{\bsnm{van~de Geer},~\bfnm{Sara~A.}\binits{S.~A.}} \AND
\bauthor{\bsnm{B{\"u}hlmann},~\bfnm{Peter}\binits{P.}}
(\byear{2009}).
\btitle{On the conditions used to prove oracle results for the {L}asso}.
\bjournal{Electron. J. Stat.}
\bvolume{3}
\bpages{1360--1392}.
\bid{doi={10.1214/09-EJS506}, issn={1935-7524}, mr={2576316}}
\end{barticle}
\bptok{imsref}%
% NOT OUTPUTED:
% issn = 1935-7524
% url = http://dx.doi.org/10.1214/09-EJS506
% fjournal = Electronic Journal of Statistics
\endbibitem


%%b15 #&#
\bibitem[\protect\citeauthoryear{van~de Geer  et al.}{2013}]{vdgetal13}
\begin{bmisc}[auto]
\bauthor{\bsnm{van~de Geer},~\bfnm{Sara}\binits{S.}}
\bauthor{\bsnm{B\"{u}hlmann},~\bfnm{P.}\binits{P.}}
\bauthor{\bsnm{Ritov},~\bfnm{Y.}\binits{Y.}}
\AND
\bauthor{\bsnm{Dezeure},~\bfnm{R.}\binits{R.}}
(\byear{2013}).
\bhowpublished{On asymptotically
optimal confidence regions and tests for high-dimensional models.
\textit{Ann. Statist.} To appear.
Available at \arxivurl{arXiv:1303.0518v2}.}
%(\byear{2013}).
\end{bmisc}
\bptok{imsref}%
%% NOT OUTPUTED:
%% issn = 0090-5364
%% url = http://dx.doi.org/10.1214/13-AOS1085
%% number = 2
%% fjournal = The Annals of Statistics
\endbibitem

%b16 #&#
\bibitem[\protect\citeauthoryear{Wasserman and Roeder}{2009}]{WR08}
\begin{barticle}[mr]
\bauthor{\bsnm{Wasserman},~\bfnm{Larry}\binits{L.}} \AND
\bauthor{\bsnm{Roeder},~\bfnm{Kathryn}\binits{K.}}
(\byear{2009}).
\btitle{High-dimensional variable selection}.
\bjournal{Ann. Statist.}
\bvolume{37}
\bpages{2178--2201}.
\bid{doi={10.1214/08-AOS646}, issn={0090-5364}, mr={2543689}}
\end{barticle}
\bptok{imsref}%
% NOT OUTPUTED:
% issn = 0090-5364
% url = http://dx.doi.org/10.1214/08-AOS646
% number = 5A
% coden = ASTSC7
% fjournal = The Annals of Statistics
\endbibitem

%b17 #&#
\bibitem[\protect\citeauthoryear{Zhang and Zhang}{2014}]{zhangzhang11}
\begin{bmisc}[auto:STB|2014/02/12|14:17:21]
\bauthor{\bsnm{Zhang},~\bfnm{C.-H.}\binits{C.-H.}} \AND
\bauthor{\bsnm{Zhang},~\bfnm{S.}\binits{S.}}
(\byear{2014}).
\bhowpublished{Confidence intervals for low-dimensional parameters with high-dimensional data.
\textit{J. R. Stat. Soc. Ser. B}.
\textbf{76}
217--242.}
\end{bmisc}
\bptok{imsref}%
\endbibitem

%b18 #&#
\bibitem[\protect\citeauthoryear{Zhao and Yu}{2006}]{zhaoyu06}
\begin{barticle}[mr]
\bauthor{\bsnm{Zhao},~\bfnm{Peng}\binits{P.}} \AND
\bauthor{\bsnm{Yu},~\bfnm{Bin}\binits{B.}}
(\byear{2006}).
\btitle{On model selection consistency of {L}asso}.
\bjournal{J. Mach. Learn. Res.}
\bvolume{7}
\bpages{2541--2563}.
\bid{issn={1532-4435}, mr={2274449}}
\end{barticle}
\bptok{imsref}%
% NOT OUTPUTED:
% issn = 1532-4435
% fjournal = Journal of Machine Learning Research (JMLR)
\endbibitem

\end{thebibliography}
\end{document}